\documentclass[11pt]{article} 
\usepackage{graphicx,type1cm,eso-pic,color}
\usepackage{amsfonts, amssymb}
\makeindex

\newcommand{\sgn}{sgn}
 
\newtheorem{theorem}{Theorem}[section]
\newtheorem{thm}{Theorem}[section]
 
\newtheorem{definition}[theorem]{Definition} 
\newtheorem{corollary}[theorem]{Corollary}

\newtheorem{cnj}[theorem]{Conjecture} 

\title{\bf Semilinear Hyperbolic Equations in   Curved Spacetime}
\author{Karen Yagdjian \\
{}\\
\small Department of Mathematics,
\small University of Texas-Pan American,\\
\small 1201 W.~University Drive, Edinburg, TX 78539,
\small USA,  
\small {yagdjian@utpa.edu}}

\begin{document}
\date{}
\maketitle

\thispagestyle{empty} 
\vspace{-0.3cm}

\begin{abstract}
This is a survey  of the author's recent work rather
than a broad survey of the literature. The survey is concerned with the global in time  solutions of the   Cauchy problem for matter waves propagating in the curved spacetimes, which can be, in particular, modeled by  cosmological models.  We   examine  the global in time  solutions of some class of semililear  hyperbolic equations, such as the Klein-Gordon equation,
which includes  the Higgs  boson   equation     in the Minkowski spacetime,   de~Sitter spacetime,  and Einstein~\&~de~Sitter spacetime.
The crucial tool for the obtaining those results is a new approach  suggested by the author based on the integral transform with the kernel containing the hypergeometric function.
\end{abstract}

\maketitle
 
\section{Introduction}
\label{S1}

\setcounter{equation}{0}
\renewcommand{\theequation}{\thesection.\arabic{equation}}

This survey is concerned with the global in time  solutions of the   Cauchy problem for matter waves propagating in the curved spacetimes, which can be, in particular, modeled by  cosmological models.
We are motivated by the significant importance of
 the qualitative description   of the global solutions of the partial differential
equations arising in the cosmological problems for  understanding of the structure of the universe and fundamental particles physics.
On the other hand, the physical implications of the mathematical results  given here  are out of the scope of this paper.
More precisely,
in this survey we   examine  the global in time  solutions of some
class of semililear  hyperbolic equations, and, in particular, the Klein-Gordon equation,
which includes  the Higgs  boson   equation     in the Minkowski spacetime,   de~Sitter spacetime,  and Einstein~\&~de~Sitter spacetime.
The Higgs  boson  plays a fundamental role in unified theories of  weak, strong, and electromagnetic
interactions \cite{Weinberg}.
\medskip

The Klein-Gordon equation arising in relativistic physics and, in particular, general relativity and cosmology, as well as, in more recent  quantum field theories,
is a covariant equation that is considered in the curved pseudo-Riemannian
manifolds.  (See, e.g., Birrell and Davies~\cite{Birrell-Davies}, Parker and Toms~\cite{{Parker-Toms}}, Weinberg~\cite{Weinberg}.)
The latest astronomical  observational discovery that the expansion of the universe is speeding   supports the model of the expanding universe that
is mathematically described by the manifold with a metric tensor depending on time and spatial variables.
\medskip

The homogeneous and isotropic cosmological models possess highest symmetry, which  makes them  more amenable to rigorous study.
Among them,  FLRW (Friedmann-Lema$\hat {\mbox{i}}$tre-Robertson-Walker) models are mentioned, which have the  flat metric of the slices of constant time.
The  FLRW spacetime metric can be written in the form
$
ds^2= -dt^2+ a^2(t)( d{x}^2   + d{y} ^2 +  d{z}^2 )\,,
$
with an appropriate scale factor $a(t)$. (See,e.g.,\cite{Hawking,Peebles,Tolman}.)
 In particular, the  metric in
de Sitter   spacetime in the Lama{\^i}tre-Robertson coordinates \cite{Moller},\cite{Tolman} has this form with the cosmic scale factor $a (t)=e^{t} $.
The time dependence of the function $a(t)$ is  determined by  the Einstein field equations for gravity  with the {\it cosmological constant} $\Lambda $,
\[
R_{\mu \nu }-\frac{1}{2} g_{\mu \nu }R = 8\pi GT_{\mu \nu }-\Lambda g_{\mu \nu }\,.
\]
The unknown of principal importance in the Einstein equations is
a metric $g$. It comprises the basic geometrical feature of the gravitational field,
and consequently explains the phenomenon of the mutual gravitational attraction of substance.
\medskip

The metric of the Einstein~\&~de~Sitter universe (EdeS universe) is a particular member of the
FLRW metrics
\begin{eqnarray*}
ds^2= -dt^2+ a^2(t)\left[ \frac{dr^2}{1-Kr^2}   + r^2  d\Omega ^2 \right]\,,
\end{eqnarray*}
where $K=-1,0$, or $+1$, for a hyperbolic, flat or spherical spatial geometry, respectively. The   Einstein~\&~de~Sitter model  of the universe is the simplest non-empty expanding model with
the line-element
\[
ds^2 = - dt^2 + a_0^2t^{4/3} \left( dx ^2+ dy ^2+dz^2 \right)
\]
in comoving coordinates \cite{Ellis}.  It was first proposed jointly by Einstein~\&~de~Sitter (the EdeS model) \cite{Einstein-Sitter}.
The observations of the microwave radiation fit in with this model  \cite{Dirac}.
The result of this case also correctly describes the early epoch, even in a universe with curvature different from zero \cite[Sec.~8.2]{Cheng}.
Even though the
EdeS spacetime is conformally flat, its causal structure
is quite different from asymptotically flat geometries. In
particular, and unlike Minkowski or Schwarzschild spacetimes,
the past particle
horizons exist.
The EdeS spacetime is a good approximation to the
large scale structure of the universe during the matter
dominated phase, when the averaged (over space and time) energy density evolves adiabatically and pressures
are vanishingly small, as, e.g., immediately after inflation. This justifies why   such a  metric  is adopted to
model the collapse of overdensity perturbations in the
early matter dominated phase that followed inflation.
\medskip

The matter waves in the  spacetime are described by the function  $\phi $, which  satisfies equations of
motion.
In the    model of universe with  curved spacetime the equation for the scalar field with  potential function \, $V$   \,
is the covariant wave equation
\[
\square_g \phi    = V'(\phi ) \quad \mbox{\rm or} \quad \frac{1}{\sqrt{|g|}}\frac{\partial }{\partial x^i}
\left( \sqrt{|g|} g^{ik} \frac{\partial \phi  }{\partial x^k} \right)  =V'(\phi ) \,,
\]
with the usual summation convention. Written explicitly in the Lama{\^i}tre- \\ Robertson coordinates  in the  de Sitter spacetime it, in particular,
for
\[
V'(\phi )=  -\mu^2  \phi + \lambda  |\phi |^{p-1}\phi,\quad   p>1,
\]
 has  the form
\begin{equation}
\label{1.1}
  \phi_{tt} +   n   \phi_t - e^{-2 t} \Delta  \phi =    \mu^2  \phi  -\lambda  |\phi |^{p-1}\phi\,,
\end{equation}
where $ \mu >0$ and $\lambda >0 $. Here $\bigtriangleup  $ is the Laplace operator on the flat metric,
$\bigtriangleup := \sum_{i=1}^n  \frac{\partial^2 }{\partial x_i^2}  $.
The equation for the Higgs real-valued scalar    field  in the
  de~Sitter spacetime is a special case of (\ref{1.1}) when $p=3$, $n=3$:
\[
  \phi_{tt} +   3   \phi_t - e^{-2 t} \Delta  \phi =    \mu^2  \phi  -\lambda   \phi^ 3  \,.
\]
Scalar fields play a fundamental role in the standard
model of particle physics, as well as its possible extensions.
In particular, scalar fields generate spontaneous
symmetry breaking and provide masses to gauge bosons
and chiral fermions by the Brout-Englert-Higgs mechanism~\cite{Englert-Brout}
using a Higgs-type potential~\cite{Higgs}.
\medskip

In the spacetime with the constant metric tensor $g$ the differential operator in the equation contains only the second-order derivatives.  For the  equation
\[
\phi_{tt}  - \Delta  \phi =    \mu^2  \phi  -\lambda  |\phi |^{p-1}\phi
\]
 the existence of a weak global solution in the energy space is known (see, e.g.,
 Proposition~3.2~\cite{G-V1989}) under certain conditions.
The equation
\begin{equation}
\label{Higgs_eq_Min}
  \phi_{tt}      -   \Delta  \phi =    \mu^2  \phi  -\lambda   \phi^ 3
\end{equation}
for the Higgs scalar    field  in the Minkowski spacetime has the time-independent flat solution
$
\phi_N (x) = \frac{\mu }{\sqrt{\lambda} } \tanh \left(   \frac{\mu  }{\sqrt{2} } N\cdot (x-x_0) \right)$, \,$ N, x_0, x \in {\mathbb R}^3$.
The unit vector $N $ defines the direction of  the propagation of the wave front.
This  solution, after Lorentz transformation, gives rise
to  a traveling solitary wave of the form
\[
\phi_{N,v} (x,t) =  \frac{\mu }{\sqrt{\lambda} } \tanh \left(   \frac{\mu  }{\sqrt{2} }  [ N \cdot (x-x_0) \pm v(t-t_0)]\frac{1}{\sqrt{1-v^2}}\right),
\,  N, x_0, x \in {\mathbb R}^3,
\]
$t \geq t_0$, if $0<v<1$, where $v $ is the initial velocity. The set of zeros of the solitary wave $\phi =\phi_{N,v} (x,t) $, that is, the set given by $N \cdot (x-x_0) \pm v(t-t_0)=0 $, is the moving boundary of the {\it wall}.  The existence of standing waves $\phi =\exp(i\omega t)v(x) $,
which are exponentially small at infinity $|x| =\infty$, and of  corresponding  solitary waves  for the equation (\ref{Higgs_eq_Min})
with $  \mu^2 <0$ and $  \lambda  <0$ is known (see, e.g.,  \cite{Strauss_77}).
\medskip

The covariant linear  wave equation  in the Einstein~\&~de~Sitter spacetime  written in the coordinates  is
\begin{eqnarray}
\label{WE}
\left(  \frac{\partial}{\partial t} \right)^2  \psi
-   t^{-{4}/{3}}  \sum_{i=1,2,3}\left( \frac{\partial }{\partial x^i} \right)^2   \psi +   \frac{2 }{ t }     \frac{\partial}{\partial t}    \psi =f\,.
\end{eqnarray}
In this survey we investigate the initial value problem for this equation and give the representation formulas for the fundamental solutions  in the case of   arbitrary dimension $n \in {\mathbb N}$
of the spatial variable $x \in {\mathbb R}^n$.
The  equation \ref{WE}    is   strictly hyperbolic in the domain with $t>0$. On the surface  $t=0$ its coefficients  have   singularities
that make   the study  of the initial value problem difficult. Then, the speed of propagation is   $t^{-\frac{2}{3}} $    for every $ t \in {\mathbb  R}\setminus \{0\}$.
The classical works on the Tricomi  and  Gellerstedt  equations
(see, e.g., \cite{Carroll-Showalter,Delache-Leray,Diaz-Weinberger,Weinstein}) appeal
to the singular Cauchy problem for the Euler-Poisson-Darboux equation,
and to the Asgeirsson mean value theorem
when handling a high-dimensional case.

\section{Method of Investigation. Integral Transform}
\setcounter{equation}{0}
\renewcommand{\theequation}{\thesection.\arabic{equation}}

We    suggested  in \cite{YagTricomi} a novel  approach to study second order hyperbolic equations with variable coefficients.
That approach avoids explicit appeal to the   Fourier integral operators, and it seems to be more immediate than the one  that uses the Euler-Poisson-Darboux equation.
It is used in a series of papers
\cite{YagTricomi}-\cite{CPDE2012}, \cite{Galstian-Kinoshita-Yagdjian}
to investigate in a unified way several equations such as the linear
and semilinear Tricomi and Tricomi-type equations, Gellerstedt
equation, the wave equation in EdeS  spacetime, the
wave and the Klein-Gordon equations in the de~Sitter and
anti-de~Sitter spacetimes.  The listed  equations play an important
role in the gas dynamics, elementary particle physics, quantum field
theory in curved spaces, and cosmology.
 For all above mentioned equations, we have obtained among  other things, fundamental solutions,   representation formulas for the initial-value problem,   $L_p-L_q$-estimates, local and global solutions for the semilinear equations, blow up phenomena, sign-changing phenomena,
self-similar solutions and number of other results.
\medskip

More precisely, in that method
the solution $v=v(x,t;b)$ to the Cauchy problem
\begin{eqnarray}
\label{1.7new}\label{1.8new}
v_{tt}-  \bigtriangleup v =0, \qquad   (t,x) \in {\mathbb R}^{1+n},  \,\,
v(x,0;b)= \varphi  (x,b), \,\, v_t(x,0)=0 ,\,\, x \in {\mathbb R}^n,
\end{eqnarray}
with the parameter $b \in B \subseteq {\mathbb R}$ is utilized. Denote that solution by $v_\varphi =v_\varphi (x,t;b)$.
There are well-known explicit representation formulas for the solution of the last problem. (See, e.g., \cite{Shatah}.)
In particular, if $\varphi $ is independent of the second variable $b$, then
$v_\varphi (x,t;b)$ does not depend on $b$ and we  write $v_\varphi (x,t)$.
\medskip

The starting point of that approach  \cite{YagTricomi} is the Duhamel's principle, which we revise in order to prepare the ground for generalization.
Our {\it first observation} is that we obtain the following representation
\begin{equation}
\label{main}
u(x,t)= \int_{t_0}^t \,d\tau \int_{ 0}^{  t-\tau  } w_f(x,z;\tau )\,dz\,,
\end{equation}
of the solution of the Cauchy problem
$
 u_{tt}-\Delta u =f(x,t)$ in ${\mathbb R}^{n+1}$, and
$u(x,t_0)=0,\quad  u_t(x,t_0)=0   \quad \mbox{\rm in}  \,\, {\mathbb R}^{n}\,,
$
where the function $w_f=w_f(x;t;\tau ) $ is the solution of the problem  \ref{1.7new}.
This formula allows us to solve problem with the source term if we
solve the problem for the {\sl same} equation without source term
but with the first initial datum.
\medskip

The  {\it second observation} is that in (\ref{main}) the upper limit $t-\tau $  of the inner integral is generated by the  propagation phenomena  with the speed which equals to one.
In fact, that is a distance function between the points at time $t$ and $\tau $.
\medskip

Our {\it third observation}  is that the solution operator $G\,:\,f \longmapsto u $ can be regarded as a composition of two operators. The first one
\[
{\mathcal W}{\mathcal E}: \,\, f  \longmapsto  w
\]
is a Fourier Integral Operator (FIO), which is a solution operator of the Cauchy problem with the first initial datum for wave equation in the Minkowski spacetime.  The second operator
\[
{\mathcal K}:\,\,w\longmapsto   u
\]
is the integral operator given by (\ref{main}). We regard the variable $z$ in  (\ref{main}) as a ``subsidiary time''. Thus, $G= {\mathcal K}\circ {\mathcal W}{\mathcal E}$ and we arrive at the diagram:
 
\includegraphics[width=0.3\textwidth]{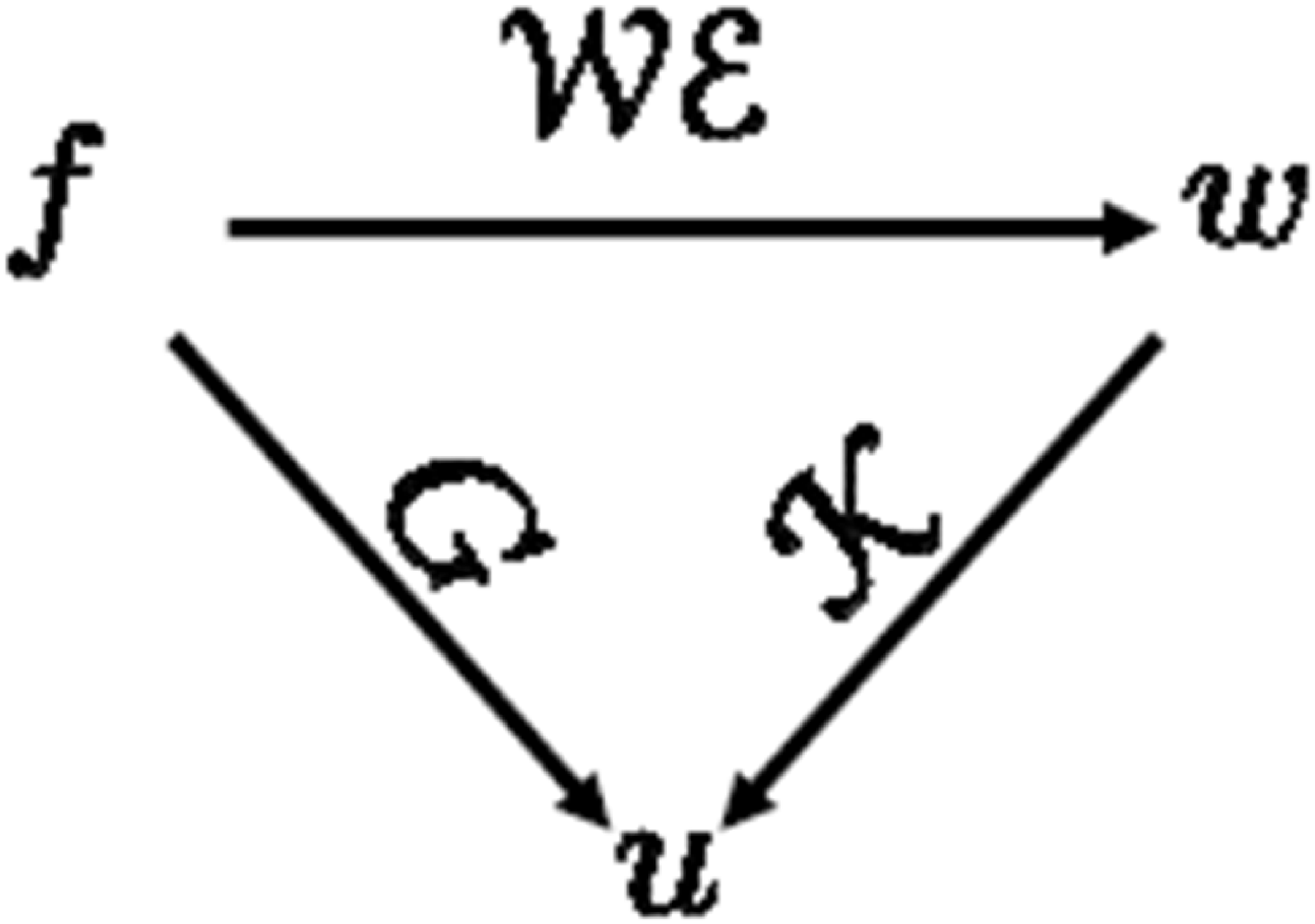}

Based on this diagram, we generated a class of operators for which we  obtained   explicit representation formulas for the solutions.
That means also that we  have representations for the {\sl fundamental solutions of the partial differential operator}. In fact,  this diagram
brings into a single  hierarchy several different partial differential operators.
Indeed, if we take into account the propagation cone by introducing the distance function
$ \phi (t)$, and if we provide   the integral operator with the kernel $ K (t;r,b) $ as follows:
\begin{equation}
\label{Oper_K}
{\mathcal K} [w](x,t)
 =
2   \int_{ t_0}^{t} db
  \int_{ 0}^{ |\phi (t)- \phi (b)|}   K (t;r,b)  w(x,r;b )  dr
, \quad x \in {\mathbb R}^n, \,\, t>t_0,
\end{equation}
then we actually can generate new representations for the solutions of  different well-known equations.
Below we give some examples of the operators with the variable coefficients. (See also \cite{yagdjian_Rend_Trieste}.)
\medskip

\paragraph{$1^0$ Tricomi-type equations.} This operator   is generated by the
kernel  $K (t;r,b)= 2E(0,t;r,b)$, where the function $  E(x,t;r,b)$
\cite{YagTricomi} is defined by
\begin{eqnarray}
\label{E}
E(x,t;r,b)
& := &
c_k \left(  (\phi (t)  + \phi (b))^2  -(x-r)^2 \right)^{-\gamma }   \\
&  &
\times F \left(\gamma , \gamma ;1; \frac{(\phi (t)  - \phi
(b))^2 - (x-r)^2} {(\phi (t)  + \phi (b))^2 - (x-r)^2}  \right), \nonumber
\end{eqnarray}
with $\gamma := k/(2k+2)$, $c_k = (k+1)^{-k/(k+1)}2^{-1/(k+1)}$, $k\not=-1$, $k  \in {\mathbb R} $,
and the distance function  is $
\phi (t)= t^{k+1}/(k+1),
$
\, while $F\big(a, b;c; \zeta \big) $ is the Gauss's hypergeometric function.
  It is proved in \cite{YagTricomi} that
for the smooth function $f=f(x,t)$, the function
\begin{eqnarray*}
u(x,t)
& = &
   2c_k \int_{ 0}^{t} db
  \int_{ 0}^{ \phi (t)- \phi (b)}  \left(  (\phi (t)  + \phi (b))^2  -r^2 \right)^{-\gamma }  \\
  &  &
\times   F \left(\gamma , \gamma ;1; \frac{(\phi (t)  - \phi
(b))^2 - r^2} {(\phi (t)  + \phi (b))^2 - r^2}  \right)  w(x,r;b )
dr, \quad  t>0,
\end{eqnarray*}
 solves the Tricomi-type equation ($l=2k \in {\mathbb N}$) (for the Tricomi equation $l=1$)
\begin{eqnarray}
\label{Tric_eq}
 u_{tt}-t^{l}\Delta u =f(x,t) \quad \mbox{\rm in}  \quad  {\mathbb R}_+^{n+1}:=\{(x,t)\,|\,  x \in {\mathbb R}^{n},\, t>0 \},
\end{eqnarray}
and takes vanishing initial values
\begin{eqnarray}
\label{vanish_initial} u(x,0)=0,\quad  u_t(x,0)=0   \quad \mbox{\rm
in}  \,\, {\mathbb R}^{n}.
\end{eqnarray}

\paragraph{$2^0$ The wave equation in the FLRW-models: de~Sitter  spacetime.} In this  example   $K (t;r,b)= 2E(0,t;r,b)$, where
the function $  E(x,t;r,b)$ \cite{Yag_Galst_CMP} is defined by
\begin{eqnarray}
\label{E0_16}
\hspace*{-0.5cm} E(x,t;r,b)
& := &
\left(  (e^{-b}  + e^{-t})^2  -(x-r)^2 \right)^{-\frac{1}{2}} \\
&  &
\times F \left(\frac{1}{2}, \frac{1}{2};1; \frac{(e^{-t}  -
e^{-b})^2 - (x-r)^2} {(e^{-t}  + e^{-b})^2 - (x-r)^2}  \right), \nonumber
\end{eqnarray}
and $ \phi (t):= 1- e^{-t}  $. For the  simplicity, in
(\ref{E0_16})  we use the notation $x^2=x\cdot x=|x|^2$ for $ x \in
{\mathbb R}^n$. It is proved in \cite{Yag_Galst_CMP} that,  defined
by the integral transform (\ref{Oper_K}) with the kernel
(\ref{E0_16})  the function
\begin{eqnarray*}
u(x,t)
& = &
 2\int_0^t \,db \int_{ 0}^{  e^{-b}-e^{-t}  } \left(  (e^{-b}  + e^{-t})^2  -r^2 \right)^{-\frac{1}{2}} \\
&  &
\times F \left(\frac{1}{2}, \frac{1}{2};1; \frac{(e^{-t}  - e^{-b})^2 - r^2}
{(e^{-t}  + e^{-b})^2 - r^2}  \right)w(x,r;\tau )\,dr
\end{eqnarray*}
 solves the wave equation in the FLRW  spaces arising in the de~Sitter
model of the universe (see, e.g. \cite{Moller}),
$ u_{tt}-e^{-2t}\Delta u =f(x,t)$   in ${\mathbb R}_+^{n+1}$,
and takes vanishing initial data (\ref{vanish_initial}).

\paragraph{$3^0$  The wave equation in the FLRW-models: anti-de~Sitter  spacetime} The third  example we obtain if we set $K (t;r,b)= 2E(0,t;r,b)$, where
the function $  E(x,t;r,b)$ is defined by (see \cite{Yag_Galst_JMAA})
\begin{eqnarray}
\label{E_017}
\hspace*{-0.5cm} E(x,t;r,b)
& := &
\left(  (e^{b}  + e^t)^2  -(x-r)^2 \right)^{-\frac{1}{2}} 
  F \left(\frac{1}{2}, \frac{1}{2};1; \frac{(e^t  -
e^{b})^2 - (x-r)^2} {(e^t  + e^{b})^2 - (x-r)^2}  \right), 
\end{eqnarray}
while the distance function is $\phi (t):= e^t-1$. In that case the
function $u=u(x,t)$  produced by the integral transform
(\ref{Oper_K}) with   $t_0=0$  and the kernel (\ref{E_017}),
solves the wave equation in the FLRW  space arising in
the anti-de~Sitter model of the universe (see, e.g. \cite{Moller}),
$
 u_{tt}-e^{2t}\Delta u =f(x,t) \quad \mbox{\rm in}  \,\, {\mathbb R}_+^{n+1}.
$
Moreover, it takes vanishing initial values (\ref{vanish_initial}).

\paragraph{$4^0$  The wave equation in the Einstein~\&~de~Sitter spacetime}
If we allow negative $k \in{\mathbb R}$ in (\ref{E}), then we obtain another way to get new operators of the above described hierarchy.
In fact,    in the hierarchy of the  hypergeometric functions $ F \left(a, b;c; \zeta  \right)$
the  simplest non-constant function is
$
  F \left(-1, -1;1; \zeta  \right) = 1+ \zeta $.
The  exponent \,$   l$\, leading to  $F \left(-1, -1;1; \zeta  \right) $   is exactly the exponent\, $   l=-4/3$ \, of the
{\sl  wave equation }(and of the {\sl  metric tensor}) {\sl  in the Einstein~\&~de~Sitter spacetime}. In that case the kernel   is
$
K(t;r,b)=
\frac{1}{18} \left(   9t^{2/3}   + 9b^{2/3}    -r^2 \right)$.
Consequently,  the function
\begin{equation}
\label{sol_EdS}
u(x,t)
=
   \int_{ 0}^{t} db
  \int_{ 0}^{ 3t^{1/3}- 3b^{1/3}}   \frac{1}{18} \left(   (3t^{1/3})^2   + (3b^{1/3})^2     -r^2 \right)  w(x,r;b )  dr
,
\end{equation}
$x \in {\mathbb R}^n$, \,$t>0$, \, solves (see \cite{Galstian-Kinoshita-Yagdjian}) the equation
 \begin{eqnarray}
 \label{0.8}
u _{tt} -  t^{-{4}/{3}}  \Delta   u  =f\quad \mbox{\rm in} \quad
{\mathbb R}_+^{n+1},
\end{eqnarray}
and  takes vanishing  initial data  (\ref{vanish_initial}) provided
that  $w ={\mathcal
WE}(f)$. Because of the singularity in the coefficient of equation
(\ref{0.8}), the Cauchy problem is not well-posed.  In order to
obtain a  well-posed problem the initial conditions must be modified
to the weighted initial value conditions
\[
\cases{
 \displaystyle   \lim_{t\rightarrow 0}\, u (x,t) = \varphi_0 (x), \qquad  \,\, x \in {\mathbb R}^n, \cr
\displaystyle
\lim_{t\rightarrow 0}
  \left(  u _t (x,t)+3 t^{ - {1}/{3}} \bigtriangleup  \varphi_0  (x   )  \right)
=  \varphi_1 (x),   \quad x \in {\mathbb  R}^n\,.
 }
\]
The operator of equation (\ref{0.8}) coincides  with the principal part of (\ref{WE}).
We remind that the   EdeS model  of
the universe
is the simplest non-empty expanding model (see, e.g., Section~4.3\cite{Ellis}).
The last equation belongs to the family of
the non-Fuchsian  partial differential equations.  There is very advanced theory of such equations (see, e.g.,
\cite{Mandai}), but according to  our knowledge the weighted initial value problem suggested in \cite{Galstian-Kinoshita-Yagdjian}
is the new one.
\medskip

More examples on the integral transforms and representation formulas are given below.

\section{Huygens' Principle  for the Klein-Gordon equation in the de~Sitter spacetime}
\setcounter{equation}{0}
\renewcommand{\theequation}{\thesection.\arabic{equation}}

In this section we  show that the   Klein-Gordon equation in the de Sitter spacetime, which is the curved manifold due to the cosmological constant, obeys the Huygens' principle only if the physical mass $m$
of the scalar field and  the dimension $n\geq 2$ of the spatial variable are tied by the equation $m^2=(n^2-1)/4 $. Recall (see, e.g., \cite{Gunther}) that a hyperbolic equation is said to satisfy Huygens'
principle if the solution vanishes at all points which cannot be
reached from the support of initial data by a null geodesic, that is, there is
no tail. The tails are important within cosmological context. (See, e.g., \cite{Ellis-Sciama,FARAONI-GUNZIG,Gleiser-Price}  and references therein.)
\medskip

Moreover,   we define the {\sl incomplete Huygens' principle}, which is the Huygens' principle
restricted to the   vanishing second initial datum,   and then prove that  the massless     scalar field  in the de Sitter spacetime
obeys the   incomplete Huygens' principle and does not obey  the  Huygens' principle, for the dimensions $n=1,3$, only.
\medskip

In  quantum field theory  for   the massive scalar field, the equation of motion is   the  Klein-Gordon equation generated by the metric $g$:
\[
\frac{1}{\sqrt{|g|}}\frac{\partial }{\partial x^i}\left( \sqrt{|g|} g^{ik} \frac{\partial \phi  }{\partial x^k} \right) = m^2 \phi  + V'(\phi ) \,.
\]
In physical
terms this equation describes a local self-interaction for a scalar particle.
In the  de~Sitter universe  the equation for the scalar field with mass   $m$   and potential   $V$
written out explicitly in coordinates is
\begin{equation}
\label{1.3}
  \phi_{tt} +   n \phi_t - e^{-2t} \bigtriangleup \phi + m^2\phi=   - V'(\phi )\,.
\end{equation}
For the solution $\Phi$ of the Cauchy problem for the linear Klein-Gordon equation
\begin{equation}
\label{1.10}
  \Phi_{tt} +   n   \Phi_t - e^{-2 t} \bigtriangleup \Phi + m^2\Phi=  0\,, \quad \Phi (x,0)= \varphi_0 (x)\, , \quad \Phi _t(x,0)=\varphi_1 (x)\,,
\end{equation}
the following formula is obtained in \cite{Yag_Galst_CMP}:
\begin{eqnarray}
\label{large}
\Phi (x,t)
& = &
e^{-\frac{n-1}{2}t} v_{\varphi_0}  (x, \phi (t)) \\
&  &
+ \,  e^{-\frac{n}{2}t}\int_{ 0}^{1} v_{\varphi_0}  (x, \phi (t)s)\big(2  K_0(\phi (t)s,t)+ nK_1(\phi (t)s,t)\big)\phi (t)\,  ds   \nonumber  \\
& &
+\, 2e^{-\frac{n}{2}t}\int_{0}^1   v_{\varphi _1 } (x, \phi (t) s)
  K_1(\phi (t)s,t) \phi (t)\, ds
, \quad x \in {\mathbb R}^n, \,\, t>0\,, \nonumber
\end{eqnarray}
provided that the mass $m$ is large, that is,  $m^2 \geq n^2/4 $.  Here, $\phi (t):= 1-e^{-t} $ and
the function   $v(x,t)$ is defined by (\ref{1.7new}). 
Next we proceed to the definition of  $K_0(z,t)  $ and $K_1(z,t)  $. (See also Section~3~\cite{Yag_Galst_CMP}.)
\medskip

We  introduce the following notations.
First, we define a  {\it chronological future}
$D_+ (x_0,t_0) $
and  a
{\it chronological past} $D_- (x_0,t_0) $  of the point   $(x_0 ,t_0)$,  $x_0 \in {\mathbb R}^n$, $t_0 \in {\mathbb R}$,
 as follows:
$
D_\pm (x_0,t_0)
  :=
 \{ (x,t)  \in {\mathbb R}^{n+1}  \, ; \,
|x -x_0 | \leq $ $\pm( e^{-t_0} - e^{-t })
\, \} $.
We define also the {\it characteristic conoid} (ray cone) by
$
C_\pm (x_0,t_0)
  :=
 \{ (x,t)  \in {\mathbb R}^{n+1}  \, ; \,
|x -x_0 | $ $ = \pm( e^{-t_0} - e^{-t })
\, \}$.
Then, we define for $(x_0, t_0) \in {\mathbb R}^n\times {\mathbb R}$  the function
\begin{eqnarray}
\label{EM}
E(x,t;x_0,t_0;M)
& :=  &
 4 ^{-M}  e^{ M(t_0+t) } \Big((e^{-t }+e^{-t_0})^2 - (x - x_0)^2\Big)^{-\frac{1}{2}+M    } \\
 &  &
\times F\Big(\frac{1}{2}-M   ,\frac{1}{2}-M  ;1;
\frac{ ( e^{-t_0}-e^{-t })^2 -(x- x_0 )^2 }{( e^{-t_0}+e^{-t })^2 -(x- x_0 )^2 } \Big) , \nonumber \\
\label{Ea}
E(x,t;x_0,t_0)
& := &
E(x,t;x_0,t_0;-iM),
\end{eqnarray}
in    $D_+ (x_0,t_0)\cup D_- (x_0,t_0) $,  where $F\big(a, b;c; \zeta \big) $ is the hypergeometric function.
The kernels  $K_0(z,t)   $, $K_1(z,t)   $, $K_0(z,t;M) $,   and  $K_0(z,t;M) $ are defined by
\begin{eqnarray}
\label{K0M}
&  &
K_0(z,t;M) \\
&  := &
4 ^ {-M}  e^{ t M}\big((1+e^{-t })^2 - z^2\big)^{  M    } \frac{1}{ [(1-e^{ -t} )^2 -  z^2]\sqrt{(1+e^{-t } )^2 - z^2} } \nonumber \\
&   &
\times  \Bigg[  \big(  e^{-t} -1 +M(e^{ -2t} -      1 -  z^2) \big)
F \Big(\frac{1}{2}-M   ,\frac{1}{2}-M  ;1; \frac{ ( 1-e^{-t })^2 -z^2 }{( 1+e^{-t })^2 -z^2 }\Big)\nonumber  \\
&  &
\hspace{0.3cm}  +   \big( 1-e^{-2 t}+  z^2 \big)\Big( \frac{1}{2}+M\Big)
F \Big(-\frac{1}{2}-M   ,\frac{1}{2}-M  ;1; \frac{ ( 1-e^{-t })^2 -z^2 }{( 1+e^{-t })^2 -z^2 }\Big) \Bigg]\nonumber ,
\end{eqnarray}
\begin{eqnarray}
\label{K1M}
K_1(z,t;M)
& :=  &
  4 ^{-M} e^{ Mt }  \big((1+e^{-t })^2 -   z  ^2\big)^{-\frac{1}{2}+M    } \\
  &  &
  \times F\left(\frac{1}{2}-M   ,\frac{1}{2}-M  ;1;
\frac{ ( 1-e^{-t })^2 -z^2 }{( 1+e^{-t })^2 -z^2 } \right). \nonumber
\end{eqnarray}
For (\ref{large}) we set  $M=\sqrt{m^2-\frac{n^2}{4}} $.
\medskip

For the case of small mass, $m^2 \leq n^2/4$, a similar formula is obtained in \cite{yagdjian_DCDS}. More precisely,
if we denote  $M=\sqrt{\frac{n^2}{4}-m^2} $, and define
\begin{eqnarray}
\label{K0}
K_0(z,t)
&   := &
K_0(z,t;-iM) ,\\
\label{K1}
K_1(z,t)
&   :=  &
 K_1(z,t;-iM)  .
\end{eqnarray}
Then for the solution $\Phi $ of the Cauchy problem  (\ref{1.10}), there is a representation
\begin{eqnarray}
\label{24}
\Phi  (x,t)
& = &
e^{-\frac{n-1}{2}t} v_{\varphi_0}  (x, \phi (t)) \\
&  &
+ \,  e^{-\frac{n}{2}t}\int_{ 0}^{1} v_{\varphi_0}  (x, \phi (t)s)\big(2  K_0(\phi (t)s,t;M)+ nK_1(\phi (t)s,t;M)\big)\phi (t)\,  ds  \nonumber \\
& &
+\, 2e^{-\frac{n}{2}t}\int_{0}^1   v_{\varphi _1 } (x, \phi (t) s)
  K_1(\phi (t)s,t;M) \phi (t)\, ds
, \quad x \in {\mathbb R}^n, \,\, t>0\,.   \nonumber
\end{eqnarray}
According to \cite{yagdjian_Rend_Trieste}, the fundamental solutions (the retarded and advanced  Green functions) of the operator have similar representations.
\medskip

Suppose now that we are looking for the simplest possible kernels $K_0(z,t;M)   $    and
$K_1(z,t;M) $   of the   integral transforms.
Surprisingly that perspective shades a light on the quantum field theory in the de~Sitter universe   and
 reveals a new unexpected  link between   the
Higuchi bound \cite{Higuchi} and the Huygens' principle.
\medskip

 Indeed, in the hierarchy of the  hypergeometric functions
the simplest one is the constant, $F \left(0, 0;1; \zeta  \right) =1$.  The parameter $M$ leading to  such function  $F \left(0, 0;1; \zeta  \right) = 1  $ is $M=\frac{1}{2}$, and, consequently, $m^2= \frac{n^2-1}{4} $.
\medskip

The next simple non-constant  function of that  hierarchy is   $F \left(-1, -1;1; \zeta  \right) = 1+ \zeta $. The   parameter $M$ leading to  such function  is     $M=\frac{3}{2} $,
and, consequently, $m^2= \frac{n^2-9}{4} $.
\medskip

In the case of $n=3$ the only real masses, which simplify the kernels, that is, make $F $ polynomial,    are      $m =\sqrt{2}$ and  $m= 0$.  These are exactly the endpoints of the
interval $(0,\sqrt{2}) $ that,  in the case of $n=3 $, is known in the quantum field theory as  the so-called Higuchi bound \cite{Higuchi}.
In fact,   the interval $  ( 0, \sqrt{2}  ) $     plays a significant role   in the linear quantum field theory ~\cite{Higuchi}, in a completely different context than the explicit representation of the solutions of the Cauchy problem.
More precisely,   the   Higuchi bound~\cite{Higuchi,Lasma Alberte,Berkhahn-Dietrichb-Hofmann,Deser-Waldron,Dengiz-Tekin},
   arises  in the quantization of free massive fields with the spin-2 in the de~Sitter spacetime with $n=3$.
It is the forbidden mass range for spin-$2$ field theory in de~Sitter spacetime because of the appearance of negative norm states.
Thus, the point $m= \sqrt{2} $ is exceptional for the   quantum fields theory in the de~Sitter spacetime.  In particular, for massive spin-2 fields, it is known \cite{Deser-Waldron,Higuchi}  that the norm of the helicity zero mode changes sign across the line $m^2= 2 $. The region  $m^2<2$ is therefore unitarily forbidden. It is noted in \cite{Lasma Alberte} that all canonically normalized helicity $-0,\pm 1,\pm 2$ modes of massive graviton on the de Sitter universe satisfy the Klein-Gordon equation for a   massive scalar field with the same effective mass.
\medskip

In the case of   $n\in {\mathbb N}$ we obtain for the physical mass several points,  $m^2= \frac{n^2}{4}-\left( \frac{1}{2} +k\right)^2$, $k=0,1,\ldots,\left[\frac{n-1}{2}\right] $,
which make $F$ polynomial.
We will call these points {\it the knot points}. For $n=1$  only  the massless field $m=0$ has a knot point.
\medskip

We state below that  the  largest knot point, and, in particular, {\sl the right endpoint of the Higuchi bound if $n=3$, is the only value of the mass of the particle which produces a
scalar field that obeys the Huygens' principle.}
\medskip

An exemplar equation satisfying Huygens' principle is the
wave equation in $n + 1$ dimensional Minkowski spacetime for odd $n
\geq  3$. According to Hadamard's conjecture (see, e.g., \cite{Gunther,Berest,{Ibragimov-Oganesyan}}) this is the only
(modulo transformations of coordinates and unknown function)
Huygensian linear second-order hyperbolic equation. There exists an extensive literature on the Huygens' principle in the 4-dimensional  spacetime of constant curvature (see e.g. \cite{FARAONI-GUNZIG,Sonego-Faraoni}  and references therein).
In \cite{ArXiV-2012}  is suggested a new proof of the following statement.
\begin{thm} \mbox{\rm \cite{ArXiV-2012}}
\label{T3}
The value $m= \sqrt{n^2-1}/2$ is the only value of the physical mass $m $, such that the solutions of the equation
\begin{eqnarray}
\label{25}
  \Phi _{tt} +   n   \Phi _t - e^{-2 t} \bigtriangleup \Phi  + m^2\Phi = 0,
\end{eqnarray}
 obey the  Huygens' principle,
whenever  the wave equation in the  Minkowski spacetime does, that is, $n \geq 3 $ is an odd number.
\end{thm}
In fact, for $m= \sqrt{n^2-1}/2$ the equation by means of the transformations of coordinates and the unknown
function can be reduced to the wave equation on the Minkowski spacetime. It is easily seen that the value  $m= \sqrt{n^2-1}/2$
is also the conformal value, i.e., the value under which a conformal change of the metric turns the problem into one on a compact in time cylinder.
\medskip

Even if the equation is not Huygensian (not tail-free for some admissible data), one might nevertheless be interested in data that produce tail-free solution.
\begin{definition} \mbox{\rm \cite{ArXiV-2012}}
We say that the equation obeys the  incomplete   Huygens' principle with respect to the first initial datum,
   if the solution with the second datum  $\varphi _1 =0 $ vanishes at all points which cannot be
reached from the support of initial data by a null geodesic.
\end{definition}
\begin{thm} \mbox{\rm \cite{ArXiV-2012}}
\label{TIHP}
Suppose that equation (\ref{25}) does not obey the   Huygens' principle. Then,  it obeys the  incomplete   Huygens' principle with respect to the first initial datum,  if and only if the equation
is massless, $m=0$, and  either $n=1 $  or $n= 3$.
\end{thm}

By combining Theorem~\ref{T3} and Theorem~\ref{TIHP} we arrive at the following interesting conclusion.

\begin{corollary} \mbox{\rm \cite{ArXiV-2012}}
Assume that the equations $
  \Phi _{tt} +   n   \Phi _t - c_1^2 e^{-2 t} \bigtriangleup \Phi  + m_1^2\Phi = 0$
and
 $ \Phi _{tt} +   n   \Phi _t - c_2^2 e^{-2 t} \bigtriangleup \Phi  + m_2^2\Phi = 0 $, where $c_1$, $c_2$ are positive numbers,
obey     the  incomplete  Huygens' principle. Then they describe the fields with different mass, $m_1\not=m_2$, (in fact, $ \frac{\sqrt{n^2-1}}{2}$
and $ 0$) if and only if the dimension $n$    is $3$.
\end{corollary}
Thus, in the de~Sitter spacetime the existence of two different scalar fields (in fact, with $m=0$ and $m^2=(n^2-1)/4 $),
which obey  incomplete Huygens' principle,
is equivalent  to the condition $n=3$.
The  dimension  $n=3$ of the last corollary agrees with the experimental data.

\section{Global Solutions of Semilinear System of Klein-Gordon Equations in de Sitter Spacetime}
\setcounter{equation}{0}
\renewcommand{\theequation}{\thesection.\arabic{equation}}

In this section we discuss  global existence of small data solutions
of the Cauchy problem for the semilinear system of Klein-Gordon
equations in the de~Sitter spacetime. Unlike   the same problem in
the Minkowski spacetime, we have no restriction on the order of the
nonlinearity  and the structure of the nonlinear term, provided that the
spectrum of the   mass matrix of the fields, which   describes the linear interactions of the fields,
is in the positive
half-line and has no intersection with some open bounded interval.
\medskip

A large amount of work has been devoted to the Cauchy problem for
the scalar semilinear Klein-Gordon equation in the Minkowski
spacetime. The existence of global weak solutions has been obtained
by  J\"orgens~\cite{Jorgens}, Segal~\cite{Segal,Segal_MSMF},
Pecher~\cite{Pecher}, Brenner~\cite{Brenner},
Strauss~\cite{Strauss},   Ginibre and Velo~\cite{G-V,G-V1989}  for
the equation
\[
  u_{tt}  -  \Delta  u  + m^2u =   |u|^{\alpha }u.
\]
For global solvability, the exact relation between $n$ and $\alpha >0$ was finally established.
More precisely,
consider the Cauchy problem for the nonlinear Klein-Gordon equation
\[
  u_{tt}  -  \Delta  u  =  - V'(u)\,,
\]
where  $V'=V'(u) $ is a nonlinear function, a typical form of which is the sum of two powers
\[
V'(u)= \lambda _0u + \lambda  |u|^{\alpha }u
\]
with $\alpha \geq 0 $ and $\lambda \geq  0$. For this equation, a
conservation of energy   is valid. For finite energy solutions
scaling arguments suggest   the assumption $ \alpha <4/(n-1)$. In
\cite{G-V1989} the existence and
uniqueness of strong global solutions in the energy space $H_{(1)}
\oplus L^2$ are proved for arbitrary space dimension  $n$ under
assumptions on $V'$ that cover the case of a sum of powers $\lambda
|u|^{\alpha }u $ with $0\leq \alpha <4/(n-1)$,  $n\geq 2  $, and $
\lambda >0$ for the highest $\alpha $. Some of the results can be
extended to the   case $\alpha =4/(n-1) $ (see, e.g. \cite{G-V},
\cite[Sec.4]{G-V1989}).
\medskip

In this section  we consider the model of interacting fields, which can be described by the system of Klein-Gordon equations with different masses,
containing interaction via mass matrix and the semilinear term.
The model obeys the following system
\begin{eqnarray}
\label{system_intr}
  \Phi _{tt} +   n  \Phi _t - e^{-2t} \bigtriangleup \Phi  + {\bf M}\Phi =   F(\Phi  )\,.
\end{eqnarray}
Here $F$ is a vector-valued function of the vector-valued function
$\Phi  $. We assume that the mass matrix $\bf M$ is real-valued,
diagonalizable, and  it has eigenvalues $m_1^2, \ldots, m_l^2$,
$i=1,2,\ldots,l$. By the similarity transformation with the
real-valued matrix $\bf O $   the mass
matrix $\bf M$  can be diagonalized. Therefore, we use the change of the
unknown function as follows:
\[
\Psi = e^{\frac{n}{2}t}{\bf O}\Phi ,\, \quad   \Phi = e^{-\frac{n}{2}t}{\bf O}^{-1}\Psi\,.
\]
Then the system (\ref{system_intr})
takes the form
\begin{eqnarray}
\label{3}
&  &
\Psi_{tt}  - e^{-2t} \bigtriangleup \Psi  + {\mathcal M}^2\Psi=    e^{\frac{n}{2}t}{\bf O}F(e^{-\frac{n}{2}t}{\bf O}^{-1}\Psi)\,,
\end{eqnarray}
where the diagonal matrix  ${\mathcal M}$, with nonnegative real part $\Re {\mathcal M} \geq 0$, is
\[
{\mathcal M}^2:= {\bf O} {\bf M}{\bf O}^{-1} - \frac{n^2}{4} {\bf I}\,, \qquad {\bf I} \quad \mbox{\rm is the identity matrix}.
\]
The matrix ${\mathcal M}^2 $ will be called the {\it curved mass matrix} of the particles, which is also sometimes referred to as  the {\it effective mass matrix}.
It is convenient to use the diagonal matrix $M= diag(|m_i^2- \frac{n^2}{4} |^{1/2} )$.
We distinguish the following three cases: the case of {\sl large mass matrix} $\bf M$ that is ${\mathcal M}^2 \geq 0 $ ($  m_i^2 \geq  \frac{n^2}{4}$,  \,$ i=1,2,\ldots, l  $);
the case of {\sl dimensional mass matrix} $\bf M$ that is
${\mathcal M}^2 = 0 $  ($  m_i^2 = \frac{n^2}{4}$,  \,$ i=1,2,\ldots, l  $); and
the case of {\sl small mass matrix} $\bf M$ that is  ${\mathcal M}^2 < 0 $ (  $  m_i^2 < \frac{n^2}{4}$,  \,$  i=1,2,\ldots, l  $). We also call the mass matrix $\bf M$ {\it critical} if $ {\mathcal
M}^2=- \frac{1}{4}I$. They lead to three different  equations:
the Klein-Gordon equation with the {\it real curved mass matrix} ${\mathcal M}$,
\begin{eqnarray*}
&  &
\Psi_{tt}  - e^{-2t} \bigtriangleup \Psi  + M^2\Psi=    e^{\frac{n}{2}t}{\bf O}F(e^{-\frac{n}{2}t}{\bf O}^{-1}\Psi)\,;
\end{eqnarray*}
the wave equation with the {\it zero curved mass matrix}
\begin{eqnarray*}
&  &
\Psi_{tt}  - e^{-2t} \bigtriangleup \Psi   =    e^{\frac{n}{2}t}{\bf O}F(e^{-\frac{n}{2}t}{\bf O}^{-1}\Psi)\,;
\end{eqnarray*}
and the Klein-Gordon equation with  the  {\it imaginary curved mass matrix}  ${\mathcal M}$,
\begin{eqnarray*}
&  &
\Psi_{tt}  - e^{-2t} \bigtriangleup \Psi  - M^2\Psi=    e^{\frac{n}{2}t}{\bf O}F(e^{-\frac{n}{2}t}{\bf O}^{-1}\Psi)\,.
\end{eqnarray*}
\medskip

Let  $ H_{(s)}
({\mathbb R}^n) $    be the Sobolev space. We use the notation $\| \cdot \|_{H_{(s)}
({\mathbb R}^n)}  $ for both the norm of a vector valued function and for the
norm of its components. To estimate the nonlinear term $F(\Phi )$
we use the following Lipschitz condition:
\medskip

\noindent {\bf Condition ($\mathcal L$).} {\it The function $F$ is
said to be Lipschitz continuous in the space $H_{(s)} ({\mathbb R}^n)$  if there are constants \,
$\alpha \geq 0$ \, and $C$ such that
\[
 \hspace{-0.4cm} \|  F(\Phi _1 )- F(\Phi _2 ) \|_{H_{(s)} ({\mathbb R}^n)}  \leq C
\|  \Phi _1 -  \Phi _2   \|_{H_{(s)} ({\mathbb R}^n)}
\Big( \|  \Phi _1  \|^{\alpha} _{H_{(s)} ({\mathbb R}^n)}
+ \|  \Phi _2   \|^{\alpha} _{H_{(s)} ({\mathbb R}^n)} \Big)
\]
for all \,\,  $\Phi _1,\Phi _2  \in H_{(s)} ({\mathbb R}^n)$. }
\medskip

If  $ s > n/2$, then any polynomial is Lipschitz continuous in the space $H_{(s)} ({\mathbb R}^n)$. For more examples of the Lipschitz continuous in the space $H_{(s)} ({\mathbb R}^n)$ functions  with low $s$  see, for example, \cite{Pecher},  \cite{Shatah}.
Define the complete metric space
\begin{eqnarray*}
&  &
X({R,s,\gamma}) \\
&  := &
\{ \Phi \in C([0,\infty) ; H_{(s)} ({\mathbb R}^n)) \; | \;
 \parallel  \Phi  \parallel _X := \sup_{t \in [0,\infty) } e^{\gamma t}  \parallel
\Phi  (x ,t) \parallel _{H_{(s)} ({\mathbb R}^n)}
\le R \}
\end{eqnarray*}
with the metric
\begin{eqnarray*}
d(\Phi _1,\Phi _2) := \sup_{t \in [0,\infty) }  e^{\gamma t}  \parallel  \Phi _1 (x , t) - \Phi _2 (x ,t) \parallel _{H_{(s)} ({\mathbb R}^n)}.
\end{eqnarray*}
The first result of this section is the following theorem.
\begin{thm} \mbox{\rm \cite{birk}}
\label{T0.1} Assume that the nonlinear term $F(\Phi )$  is
Lipschitz continuous in the  space $H_{(s)} ({\mathbb R}^n)$, $ s > n/2\geq
1$,  $\alpha >0 $, and $F(0)=0$. Assume also that the system has a
large mass matrix. Then, there exists $\varepsilon _0>0 $ such that,
for every given vector-valued functions $\varphi_0 ,\varphi_1 \in {
H}_{(s)} ({\mathbb R}^n)$,  such that
\[
 \| \varphi_0   \|_{ { H}_{(s)} ({\mathbb R}^n)}
+  \|\varphi_1  \|_{ { H}_{(s)} ({\mathbb R}^n)} \leq \varepsilon, \qquad \varepsilon  < \varepsilon_0\,,
\]
there exists a global solution $\Phi \in C^1([0,\infty);H_{(s)} ({\mathbb R}^n))$ of the Cauchy problem
\begin{eqnarray}
\label{NWE}
&  &
  \Phi _{tt} +   n   \Phi _t - e^{-2 t} \Delta  \Phi  +  {\bf M}  \Phi  =   F(\Phi )\,,\\
\label{ICPHI}
&  &
 \Phi  (x,0)= \varphi_0 (x)\, , \quad \Phi  _t(x,0)=\varphi_1 (x) \,.
 \end{eqnarray}
The  solution  $ \Phi  (x ,t) $  belongs to the space  $  X({2\varepsilon,s,0})  $, that is,
\[
\sup_{t \in [0,\infty)} \|\Phi  (x ,t) \|_{H_{(s)} ({\mathbb R}^n)}  < 2\varepsilon .
\]
 \end{thm}

For the scalar equation   this theorem implies   Theorem~0.1~\cite{JMAA2012}.
In fact, for the scalar equation  if
$
F(\Phi )= \pm  |\Phi |^\alpha \Phi$   or $F(\Phi )= \pm  |\Phi |^{\alpha+1}$,
then, according to Theorem~0.1~\cite{JMAA2012}, the small data Cauchy problem is  globally  solvable for every $\alpha \in (0,\infty) $ if
$m \in (0, \sqrt{n^2-1}/2)\cup [n/2,\infty) $ and the condition (${\mathcal L}$)
is fulfilled.
\begin{cnj} \mbox{\rm \cite{birk}}
The open interval  $  (  \sqrt{n^2-1}/2,n/2 ) $
is a forbidden physical mass interval for the small data global solvability of the Cauchy problem for all  $\alpha \in (0,\infty) $.
\end{cnj}

Consider the particular case of the scalar equation with the spatial
dimension $n=3$. In this case the interval  $  (  \sqrt{n^2-1}/2,n/2
) $ for the physical mass  is reduced to  $  (  \sqrt{2} ,3/2 ) $.
 The interval $  ( 0, \sqrt{2}  ) $    is the   Higuchi
bound.
This is why  we pay special attention to the system of equations with the mass matrix  $\bf M $ which is
 orthogonally similar to the matrix $ \frac{n^2-1 }{4}{\bf I}$. We  call such mass matrix  $\bf M $ critical.
We also call the mass matrix  $\bf M$ {\it semi-critical mass matrix}  if
  the spectrum $\sigma ({\bf M}) $ of the mass matrix $\bf M$ is a subset of  $  (0, (n^2-1)/4 ] $.
For the system with the semi-critical  mass matrix  $\bf M $ we prove the   global existence, which is not known in the critical case even for the
scalar equation.

\begin{thm}\mbox{\rm \cite{birk}}
\label{T0.2}
Assume that the nonlinear term $F(\Phi )$  is  Lipschitz continuous in the  space $H_{(s)} ({\mathbb R}^n)$, $ s > n/2\geq 1$,   $\alpha >0 $, and $F(0)=0$.
Assume also that the mass matrix  $\bf M $ is semi-critical, that is  $\sigma ({\bf M})\subset  (0, (n^2-1)/4 ]$.
Then, there exists $\varepsilon _0>0 $ such that, for every given vector-valued functions $\varphi_0 ,\varphi_1 \in H_{(s)} ({\mathbb R}^n)$,  such that
\[
 \| \varphi_0   \|_{ { H}_{(s)} ({\mathbb R}^n)}
+  \|\varphi_1  \|_{ { H}_{(s)} ({\mathbb R}^n)} \leq \varepsilon, \,\quad \varepsilon  < \varepsilon_0\,,
\]
there exists a global solution $\Phi \in C^1([0,\infty);H_{(s)} ({\mathbb R}^n))$ of the Cauchy problem (\ref{NWE}), (\ref{ICPHI}).
The  solution \, $ \Phi  (x ,t) $ \, belongs to the space \, $  X({2\varepsilon,s,\gamma })  $,
where
\[
\gamma <  \frac{1}{\alpha +1} \left(\frac{n }{2} -\max \left\{  \sqrt{\frac{n^2}{4}-\lambda }\,; \lambda \in\sigma ({\bf M}) \right\} \right),
\]
that is,
\[
\sup_{t \in [0,\infty)}  e^{\gamma t}  \|\Phi  (x ,t) \|_{H_{(s)} ({\mathbb R}^n)}  < 2\varepsilon  .
\]
 \end{thm}

Baskin~\cite{BaskinSE} discussed small data global solutions for the
scalar Klein-Gordon equation on asymptotically de Sitter spaces,
which is a   compact manifold with boundary.   More precisely, in
\cite{BaskinSE}  the Cauchy problem is considered for the
semilinear equation
$
\Box_g u +m^2 u = f(u)$, $
u(x,t_0)=\varphi_0 (x) \in H_{(1)}({\mathbb R}^n)$ , $u_t(x,t_0)=\varphi_1 (x) \in L^2 ({\mathbb R}^n)$,
where mass is large, $m^2> n^2/4$.   In
Theorem~1.3 \cite{BaskinSE} the existence of the global solution for
small energy data is stated. (For  references on the
asymptotically de Sitter spaces, see
\cite{Baskin},
 \cite{Vasy_2010}.)

\section{The Scalar Equation.  Case of  Large   Mass}
\label{SLM}
\label{sec:2}
\setcounter{equation}{0}
\renewcommand{\theequation}{\thesection.\arabic{equation}}

We extract a linear part of the system    (\ref{3}) as an initial model that must be treated
first. That linear system is diagonal, which allows us to restrict ourselves to one scalar equation
\begin{equation}
\label{K_G_linear}
u_{tt} - e^{-2t} \bigtriangleup u  +{\mathcal M}^2 u=  f\,,
\end{equation}
where  ${\mathcal M}  $ is a non-negative number throughout this section.
The  equation (\ref{K_G_linear}) is strictly hyperbolic. That
implies the well-posedness of the Cauchy problem for
(\ref{K_G_linear}) in several function spaces. The coefficients of
the equation are   analytic functions and, consequently, Holmgren's
theorem  implies  local uniqueness in the space of distributions.
Moreover, the speed of propagation is  equal to
$e^{-t} $ for every $ t \in {\mathbb R}$. The second-order strictly
hyperbolic  equation (\ref{K_G_linear}) possesses two fundamental
solutions resolving the Cauchy problem. They can be written
 in terms of Fourier integral operators \cite{Horm},
which give a complete description of the wave front sets of the
solutions. The distance $
e^{-t}|\xi |$ between two characteristic roots  of the equation
(\ref{K_G_linear})  tends to zero as $t \to +\infty
$. Thus, the operator is not uniformly
strictly hyperbolic. The finite integrability of the
characteristic roots  leads to the existence of a so-called {\it horizon} for
that equation.  The
equation (\ref{K_G_linear}) is neither Lorentz  invariant nor
invariant with respect to  scaling and that brings additional
difficulties.
\medskip

\noindent
{\bf $\bf L^p-L^q$ estimates for  equation with  source.}
We consider the  equation with $n\geq 2$.
The solution $u= u(x,t)$ to the Cauchy problem
\begin{equation}
\label{1.26}
u_{tt} - e^{-2t}\Delta u +M^2 u= f ,\quad u(x,0)= 0  , \quad u_t(x,0)=0,
\end{equation}
with \, $ f \in C^\infty ({\mathbb R}^{n+1})$\, and with   vanishing
initial data is given by the next expression
\[
u(x,t)
   =
2   \int_{ 0}^{t} db
  \int_{ 0}^{ e^{-b}- e^{-t}} dr \, \,  v_f(x,r ;b) E(r,t; 0,b)  ,
\]
where the function
$v_f(x,t;b)$
is a solution to the Cauchy problem
\begin{equation}
\label{1.6}
v_{tt} -   \bigtriangleup v  =  0 \,, \quad v(x,0;b)=f(x,b)\,, \quad v_t(x,0;b)= 0\,.
\end{equation}
 Thus, for the solution $\Phi$ of the equation
\begin{equation}
\label{1.7a}
  \Phi_{tt} +   n  \Phi_t - e^{-2 t} \bigtriangleup \Phi + m^2\Phi=  f\,,
\end{equation}
due to the relation $u = e^{\frac{n}{2}t}\Phi$, we obtain
\begin{eqnarray}
\label{1.29l}
\Phi  (x,t)
&  =  &
2  e^{-\frac{n}{2}t} \int_{ 0}^{t} db
  \int_{ 0}^{ e^{-b}- e^{-t}} dr  \,  e^{\frac{n}{2}b}v_f(x,r ;b) E(r,t; 0,b) .
\end{eqnarray}
For the solution $u=u(x,t)$   of the Cauchy problem (\ref{1.26})  according to Corollary~9.3~\cite{Yag_Galst_CMP}
  one has estimate
\begin{eqnarray*}
&  &
\| (-\bigtriangleup )^{-s} u(x,t) \|_{ { L}^{  q} ({\mathbb R}^n)  }  \\
& \le &
C_M\int_{ 0}^{t} \|  f(x, b)  \|_{ { L}^{p} ({\mathbb R}^n)  }
 e^{b}\left( e^{-b}- e^{-t}  \right)^{ 1+  2s-n(\frac{1}{p}-\frac{1}{q}) } \left(1+  t-  b    \right)^{1- \sgn M}\,  db\,,
\end{eqnarray*}
provided that  $1<p \leq 2$, $\frac{1}{p}+ \frac{1}{q}=1$, $\frac{1}{2} (n+1)\left( \frac{1}{p} - \frac{1}{q}\right) \leq
2s \leq n \left( \frac{1}{p} - \frac{1}{q}\right)< 2s  +1 $. Thus, for the solution $\Phi $ (\ref{1.29l}) of the equation (\ref{1.7a}), we obtain
\begin{eqnarray*}
&  &
\| (-\bigtriangleup )^{-s} \Phi (x,t) \|_{ { L}^{  q} ({\mathbb R}^n)  }
 \le
C_Me^{-\frac{n}{2}t}\\
&  &
\qquad \times \int_{ 0}^{t}e^{\frac{n}{2}b} \|  f(x, b)  \|_{ { L}^{p} ({\mathbb R}^n)  }
 e^{b}\left( e^{-b}- e^{-t}  \right)^{ 1+  2s-n(\frac{1}{p}-\frac{1}{q}) } \left(1+  t-  b    \right)^{1- \sgn M}\,  db .
\end{eqnarray*}
In particular, for $s=0$ and $p=q=2$,    we have
\[
\| \Phi  (x,t) \|_{ { L}^{  2} ({\mathbb R}^n)  }
\le
C_Me^{-\frac{n}{2}t}\int_{ 0}^{t}e^{\frac{n}{2}b} \|  f(x, b)  \|_{ { L}^{2} ({\mathbb R}^n)  }   \left(1+  t-  b    \right)^{1- \sgn M}
 \,  db.
\]
Here the rates of exponential factors are independent of  ${\mathcal M}$ and, consequently, of the mass $m$.

\paragraph{$\bf L^p-L^q$ estimates for  equations without source.}
According to Theorem~10.1~\cite{Yag_Galst_CMP} the solution  $u=u(x,t)$ of the Cauchy problem
\[
u_{tt}-  e^{-2t} \bigtriangleup u +M^2 u =0\,, \quad u(x,0)= \varphi_0 (x)\, , \quad u_t(x,0)=\varphi_1 (x)\,,
\]
 satisfies the following $L^p-L^q$ estimate
\begin{eqnarray*}
\| (-\bigtriangleup )^{-s} u(x,t) \|_{ { L}^{  q} ({\mathbb R}^n)  }
& \le &
C_M  (1+ t  )^{1- \sgn M}(1- e^{-t}) ^{ 2s-n(\frac{1}{p}-\frac{1}{q}) }\\
& &
\times \Big\{ e ^{\frac{t}{2}}  \| \varphi_0   \|_{ { L}^{p} ({\mathbb R}^n)  }
+ (1- e^{-t}) \|\varphi_1  \|_{ { L}^{p} ({\mathbb R}^n)  }
 \Big\}
\end{eqnarray*}
for all $t \in (0,\infty)$, provided that  $1<p \leq 2$, $\frac{1}{p}+ \frac{1}{q}=1$, $\frac{1}{2} (n+1)\left( \frac{1}{p} - \frac{1}{q}\right) \leq
2s \leq n \left( \frac{1}{p} - \frac{1}{q}\right)  < 2s +1 $.
In particular,   for large $t $   we obtain the following {\it no decay} estimate
\[
\| (-\bigtriangleup )^{-s} u(x,t) \|_{ { L}^{  q} ({\mathbb R}^n)  }
  \leq
C_M (1+ t  )^{1- \sgn M} \Big\{ e ^{\frac{t}{2}}  \| \varphi_0   \|_{ { L}^{p} ({\mathbb R}^n)  }
+  \|\varphi_1  \|_{ { L}^{p} ({\mathbb R}^n)  }
 \Big\}\,.
\]
Thus, for the solution $\Phi $ of the Cauchy problem (\ref{1.10}),
due to the relation $u = e^{\frac{n}{2}t}\Phi $, we obtain the decay  estimate
\begin{eqnarray*}
\| (-\bigtriangleup )^{-s} \Phi (x,t) \|_{ { L}^{  q} ({\mathbb R}^n)  }
&  \leq &
C_M e^{-\frac{n}{2}t} (1+ t  )^{1- \sgn M}(1- e^{-t}) ^{ 2s-n(\frac{1}{p}-\frac{1}{q}) } \\
&  &
\times \Big\{ e ^{\frac{t}{2}}  \| \varphi_0    \|_{ { L}^{p} ({\mathbb R}^n)  }
+ (1- e^{-t}) \|\varphi_1  \|_{ { L}^{p} ({\mathbb R}^n)  }
 \Big\}
\end{eqnarray*}
for all $t>0$.

\section{The Scalar Equation. Imaginary Curved Mass}
\label{S2}
\label{sec:3}
\setcounter{equation}{0}
\renewcommand{\theequation}{\thesection.\arabic{equation}}

In this section we consider the linear part of the scalar equation
\begin{equation}
\label{K_G_Higgs}
u_{tt} - e^{-2t} \bigtriangleup u  - M^2   u=  - e^{\frac{n}{2}t}V'(e^{-\frac{n}{2}t}u ),
\end{equation}
with $M\geq 0 $. The  equation (\ref{K_G_Higgs}) covers two important cases. The first one is the Higgs boson equation, which has $V'(\phi )=\lambda \phi ^3 $
and $M^2=  \mu m^2+ n^2/4 $ with $\lambda >0 $ and $\mu >0 $, while $n=3$.
The second case is the case of the small physical  mass, that is $0 \leq m  \le \frac{n }{2}$. For the last case
$ M^2= \frac{n^2}{4}-m^2$.
The solution $u= u(x,t)$ to the Cauchy problem
\begin{eqnarray}
\label{eqim}
u_{tt} - e^{-2t}\Delta u -M^2 u= f ,\quad u(x,0)= 0  , \quad u_t(x,0)=0,
 \end{eqnarray}
with \, $ f \in C^\infty ({\mathbb R}^{n+1})$\, and with   vanishing
initial data is given in \cite{yagdjian_DCDS} by the next expression
\begin{eqnarray*}
u(x,t)
&  =  &
2   \int_{ 0}^{t} db
  \int_{ 0}^{ e^{-b}- e^{-t}} dr  \, \, v(x,r ;b) E(r,t; 0,b;M)  ,
\end{eqnarray*}
where the function
$v(x,t;b)$
is a solution to the Cauchy problem for the  wave equation (\ref{1.6}).
\medskip

The solution $u=u (x,t)$ to the Cauchy problem
\[
u_{tt}-  e^{-2t} \bigtriangleup u -M^2 u =0\,, \quad u(x,0)= \varphi_0 (x)\, , \quad u_t(x,0)=\varphi_1 (x)\,,
\]
with \, $ \varphi_0 $,  $ \varphi_1 \in C_0^\infty ({\mathbb R}^n) $, $n\geq 2$, can be represented (see \cite{yagdjian_DCDS}) as follows:
\begin{eqnarray*}
u(x,t)
& = &
 e ^{\frac{t}{2}} v_{\varphi_0}  (x, \phi (t))
+ \, 2\int_{ 0}^{1} v_{\varphi_0}  (x, \phi (t)s) K_0(\phi (t)s,t;M)\phi (t)\,  ds  \nonumber \\
& &
+\, 2\int_{0}^1   v_{\varphi _1 } (x, \phi (t) s)
  K_1(\phi (t)s,t;M) \phi (t)\, ds
, \quad x \in {\mathbb R}^n, \,\, t>0\,,
\end{eqnarray*}
where $\phi (t):=  1-e^{-t} $.
 Here, for $\varphi \in C_0^\infty ({\mathbb R}^n)$ and for $x \in {\mathbb R}^n$,
the function $v_\varphi  (x, t)$  is   the solution   of the Cauchy problem (\ref{1.7new}).
Thus, for the solution $\Phi $ of the Cauchy problem
\begin{eqnarray*}
  \Phi _{tt} +   n   \Phi _t - e^{-2 t} \bigtriangleup \Phi  + m^2\Phi =  f ,\quad \Phi  (x,0)= 0  , \quad \Phi  _t(x,0)=0,
\end{eqnarray*}
due to the relation $u = e^{\frac{n}{2}t}\Phi $, we obtain
with \, $ f \in C^\infty ({\mathbb R}^{n+1})$\, and with   vanishing
initial data   the next expression
\begin{eqnarray}
\label{Pfim}
\Phi  (x,t)
  =
2   e^{-\frac{n}{2}t}\int_{ 0}^{t} db
  \int_{ 0}^{ e^{-b}- e^{-t}} dr  \,  e^{\frac{n}{2}b} v_f(x,r ;b) E(r,t; 0,b;M),
\end{eqnarray}
where the function
$v_f(x,t;b)$
is a solution to the Cauchy problem for the  wave equation (\ref{1.6}).
Thus, for the solution $\Phi $ of the Cauchy problem  (\ref{1.10}),
due to the relation $u = e^{\frac{n}{2}t}\Phi $, we obtain (\ref{24}).
In fact, the representation formulas of this section have been used in \cite{CPDE2012} to establish sign-changing
properties of the solutions of the Higgs boson equation.

\paragraph{The critical case of $\bf m^2= (n^2-1)/4 $.}
In   this case we have
\begin{eqnarray*}
E\left(x,t;x_0,t_0;\frac{1}{2}\right) =
 \frac{1}{2} e^{ \frac{1}{2}(t_0+t) } ,\quad E\left(z,t;0,b;\frac{1}{2}\right) = \frac{1}{2} e^{ \frac{1}{2}(b+t) }\,,
\end{eqnarray*}
while
\begin{eqnarray*}
K_0\left(z,t;\frac{1}{2} \right)
  =
- \frac{1}{4}  e^{ \frac{1}{2}t }  ,\qquad
K_1\left(z,t;\frac{1}{2} \right)
   =   \frac{1}{2}  e^{ \frac{1}{2}t }   \,.
 \end{eqnarray*}
For the solution (\ref{1.29l}) of the equation (\ref{1.7a}) with the source term it follows
\begin{eqnarray*}
\Phi  (x,t)
&  =  &
  e^{-\frac{n-1}{2}t}\int_{ 0}^{t}  e^{\frac{n+1}{2}b}
   V_{f}(x,e^{-b}- e^{-t};b)   \,   db .
\end{eqnarray*}
where  $v(x,r ;b) $ is defined by (\ref{1.6}).
We denote by $V_{f}(x,t;b)$ the solution of the problem
$
V_{tt}-  \bigtriangleup V =0, \quad V(x,0)= 0, \quad V_t(x,0)=f (x,b)
$.
Further, for the solution $\Phi $ (\ref{24}) of the  equation  without source term we have
\begin{eqnarray*}
\Phi  (x,t)
 & =  &
e^{-\frac{n-1}{2}t}  \left(  \frac{\partial V_{\varphi _0 }}{\partial t} \right) (x, 1-e^{-t}) \\
&  &
+ \,   \frac{n-1}{2}e^{-\frac{n-1}{2}t}    V_{\varphi_0}  (x, 1-e^{-t} )
+\,  e^{-\frac{n-1}{2}t}  V_{\varphi_1}  (x, 1-e^{-t} ),
\end{eqnarray*}
$x \in {\mathbb R}^n, \,\, t>0$, where  we denote by $V_{\varphi }$ the
solution of the problem
$
V_{tt}-  \bigtriangleup V =0, \quad V(x,0)= 0, \quad V_t(x,0)=\varphi (x)$.
Thus, in particular, we arrive at the next theorem.

\begin{thm} \mbox{\rm \cite{birk}}
\label{CHP}
The  solutions of the equation
$
  \Phi _{tt} +   n   \Phi _t - e^{-2 t} \Delta  \Phi  +  {\bf M}  \Phi  =  0\,,
 $
 obey the strong Huygens' Principle, if and only if $n \geq 3 $ is an odd number and the mass matrix $\bf M $ is  the diagonal
 matrix $\frac{n^2-1}{4}{\bf I} $.
\end{thm}

\begin{thm}\mbox{\rm \cite{yagdjian_DCDS}}
\label{L1} Suppose that  $m^2= ( n^2-1 )/4$. If $ \varphi _0
=\varphi _1 =0$ and $\frac{1}{2}(n+1)(\frac{1}{p}-\frac{1}{q}) -1
\le 2s \le n (\frac{1}{p}-\frac{1}{q})$, then for the solution $\Phi
= \Phi (x,t)$ of the equation (\ref{1.7a}) the following estimate
holds
\begin{eqnarray*}
&  &
\|
(-\bigtriangleup )^{-s} \Phi  (x,t) \|_{ { L}^{q} ({\mathbb R}^n)  }  \\
&  \leq  &
C e^{-\frac{n-1}{2}t}\int_{ 0}^{t}  e^{\frac{n+1}{2}b}
(e^{-b}- e^{-t})^{1+2s-n(\frac{1}{p}-\frac{1}{q})}  \|
  f(x, b)   \|_{ { L}^{p} ({\mathbb R}^n)  }    db
, \quad  t>0\,.
\end{eqnarray*}
For the solution $\Phi = \Phi (x,t)$  of the  Cauchy problem
(\ref{1.10}): if $ f \equiv 0$, $ \varphi _0 =0$, and
$\frac{1}{2}(n+1)(\frac{1}{p}-\frac{1}{q}) -1 \le 2s \le n
(\frac{1}{p}-\frac{1}{q})$, then
\begin{eqnarray*}
\|
(-\bigtriangleup )^{-s} \Phi  (x,t) \|_{ { L}^{q} ({\mathbb R}^n)  }
& \leq  &
C e^{-\frac{n-1}{2}t} (1-e^{-t})^{1+2s-n(\frac{1}{p}-\frac{1}{q})}  \| \varphi _1   \|_{ { L}^{p}({\mathbb R}^n)  }
, \quad  t>0\,,
\end{eqnarray*}
while if $ f \equiv 0$, $ \varphi _1 =0$, and
$\frac{1}{2}(n+1)(\frac{1}{p}-\frac{1}{q})  \le 2s \le n
(\frac{1}{p}-\frac{1}{q})$, then
\begin{eqnarray*}
\|
(-\bigtriangleup )^{-s} \Phi  (x,t) \|_{ { L}^{q} ({\mathbb R}^n)  }
& \leq  &
C e^{-\frac{n-1}{2}t} (1-e^{-t})^{2s-n(\frac{1}{p}-\frac{1}{q})} \| \varphi _0   \|_{ { L}^{p}({\mathbb R}^n)  }
, \quad  t>0\,.
\end{eqnarray*}
\end{thm}
\medskip

To complete the list of the $ L^p-L^q$ estimates we quote below
results (Theorems~2.2, 2.6  from \cite{JMAA2012}), which are applicable to the scalar
equation with noncritical mass.
The solution  $\Phi =\Phi (x,t)$ of the Cauchy problem
\[
  \Phi _{tt} +   n   \Phi _t - e^{-2 t} \bigtriangleup \Phi  \pm m^2\Phi =  0\,, \quad \Phi  (x,0)= \varphi_0 (x)\, , \quad \Phi  _t(x,0)=\varphi_1 (x)\,,
\]
with either $M=\sqrt{\frac{n^2}{4}-m^2}$  and  $m < \sqrt{n^2-1}/{2}$ for the case of ``plus'', or  $M=\sqrt{\frac{n^2}{4}+m^2}$ for the case of ``minus'',
satisfies the following $L^p-L^q$ estimate
\begin{eqnarray*}
&  &
\| (-\bigtriangleup )^{-s}\Phi (x,t) \|_{ { L}^{  q} ({\mathbb R}^n)  }\\
 &  \leq &
C_{M,n,p,q,s}   (1- e^{-t}) ^{ 2s-n(\frac{1}{p}-\frac{1}{q}) }e^{( M -\frac{n}{2})t}    \Big\{   \| \varphi_0   \|_{ { L}^{p} ({\mathbb R}^n)  }
+ (1- e^{-t}) \|\varphi_1  \|_{ { L}^{p} ({\mathbb R}^n)  }
 \Big\}
\end{eqnarray*}
for all $t \in (0,\infty)$, provided that  $1<p \leq 2$, $\frac{1}{p}+ \frac{1}{q}=1$, $\frac{1}{2} (n+1)\left( \frac{1}{p} - \frac{1}{q}\right) \leq
2s \leq n \left( \frac{1}{p} - \frac{1}{q}\right)  < 2s +1 $.
\medskip

Further, let $\Phi =\Phi (x,t)$ be the solution of the Cauchy
problem
\begin{equation}
\label{2.5}
  \Phi _{tt} +   n   \Phi _t - e^{-2 t} \bigtriangleup \Phi  \pm m^2\Phi =  f\,, \quad \Phi  (x,0)= 0\, , \quad \Phi  _t(x,0)=0\,,
\end{equation}
with either $M=\sqrt{\frac{n^2}{4}-m^2} $  and  $m < \sqrt{n^2-1}/2$ for the case of ``plus'', or  $M=\sqrt{\frac{n^2}{4}+m^2} $ for the case of ``minus''.
Then  $\Phi =\Phi (x,t)$ satisfies the following $L^p-L^q$ estimate:
\begin{eqnarray*}
&  &
\| (-\bigtriangleup )^{-s}  \Phi (x ,t) \|_{ { L}^{q} ({\mathbb R}^n)  }\\
&  \leq   &
 C_M e^{ -M t  } e^{-\frac{n}{2}t}  e^{-t[2s-n(\frac{1}{p}-\frac{1}{q})]}  \\
 &  &
\times \int_{ 0}^{t} e^{\frac{n}{2}b} e^{ M b  }(e^{t-b}-1)^{1+2s-n(\frac{1}{p}-\frac{1}{q})} (e^{t-b}+1)^{2 M-1}\| f(x ,b)  \|_{ { L}^{p}({\mathbb R}^n)  }  \,db\,,
\end{eqnarray*}
as well as
\begin{eqnarray*}
\| (-\bigtriangleup )^{-s}  \Phi (x ,t) \|_{ { L}^{q} ({\mathbb R}^n)  } 
&  \leq   &
 C_M   e^{-(\frac{n}{2}-M)t}
\int_{ 0}^{t} e^{(\frac{n}{2}-M)b} e^{-b(2s-n(\frac{1}{p}-\frac{1}{q}))}  \| f(x ,b)  \|_{ { L}^{p}({\mathbb R}^n)  }  \,db 
\end{eqnarray*}
for all $t\in (0,\infty)$, and for the above written range of the parameters $p$, $q$, $s$.

\section{Global Existence. Small Data Solutions}
\label{S3}
\label{sec:4}
\setcounter{equation}{0}
\renewcommand{\theequation}{\thesection.\arabic{equation}}

The Cauchy problem (\ref{eqim}) for the scalar equation  was studied
in \cite{yagdjian_DCDS}. For  $F(\Phi) = c
|\Phi  |^{\alpha +1} $, $c\not=0$, Theorem~1.1~\cite{yagdjian_DCDS},
implies  nonexistence of a global solution even for arbitrary small
initial data $\varphi_0  (x ) $ and  $\varphi_1  (x ) $ under
some conditions on $n$, $\alpha $, and $M$.
Theorem~3.1\cite{JMAA2012} gives  the  blow
up result for the equation with imaginary
physical mass.
That theorem shows that instability of the trivial solution occurs in a very strong
sense, that is,  an arbitrarily small perturbation of the initial data can make the perturbed
solution blowing up in finite time.
If we allow large initial data, then, according to Theorem~1.2~\cite{yagdjian_DCDS}, the concentration of the mass, due to the non-dispersion property of
the  de~Sitter spacetime, leads to the nonexistence of the global solution, which
cannot be recovered even by adding an exponentially decaying factor in the nonlinear term.
\medskip

In this section we are going to study the global existence of
solutions for the system of semilinear Klein-Gordon equations. The
first step toward such result is to establish the
$L_p-L_q$-estimates for the equation with source term. For the
scalar equation this estimate is proved in \cite{JMAA2012}.
Although we want to prove the global existence for two different cases, for the system with the semi-critical mass matrix and for the system of equations with the
large mass matrix,
the consideration  can be done in the single framework.
\medskip

We reduce the Cauchy problem to the integral
equation. The main tool for such reduction is the fundamental
solution  for the interacting fields, which
can be described by the system of Klein-Gordon equations  containing
interaction via mass matrix and the semilinear term. The model obeys
the  system (\ref{system_intr}).
By the similarity transformation $\bf O $ the mass matrix $ \bf M$
can be diagonalized, therefore we use a change of unknown function
\[
\Psi ={\bf O} \Phi\,,\quad \Phi ={\bf O} ^{-1}\Psi,
\]
and arrive at
\begin{eqnarray*}
&  &
\Psi_{tt} +   n H \Psi _t - e^{-2Ht} \bigtriangleup \Psi  +  \widetilde{\bf M}\Psi=    \widetilde{F}(\Psi)\,,
\end{eqnarray*}
where
\[
\widetilde{\bf M} :={\bf O} M{\bf O} ^{-1}=  \mbox{diag}\{m_1^2,\ldots,m_l^2\} ,
 \quad  \widetilde{F}(\Psi) := {\bf O} F({\bf O} ^{-1}\Psi ) .
\]
Let us consider the linear diagonal  system
\begin{eqnarray*}
&  &
\Psi_{tt} +   n H \Psi _t - e^{-2Ht} \bigtriangleup \Psi  +  \widetilde{\bf M}\Psi=    \widetilde{f}\,.
\end{eqnarray*}
Here $\widetilde{f} $ is a vector-valued function with the components $f_i $, $i=1,\ldots,l$.
Then, the solution of the Cauchy problem for the last system with the initial conditions
\[
\Psi (x,0)=   0,\qquad \Psi_t (x,0)=   0 ,
\]
is
\begin{eqnarray*}
\Psi   (x,t)
&  =  &
2  e^{-\frac{n}{2}t} \int_{ 0}^{t} db
  \int_{ 0}^{ e^{-b}- e^{-t}} dr  \,  e^{\frac{n}{2}b}\widetilde{E} (r,t; 0,b) \widetilde{v}(x,r ;b)  ,
\end{eqnarray*}
where the components $v_i $, $i=1,\ldots,l$,  of the vector-valued function
$\widetilde{v}(x,t;b)$
are   solutions to the Cauchy problem for the  wave equation
\begin{equation}
\label{WEsystem}
v_{tt} -   \bigtriangleup v  =  0 \,, \quad v(x,0;b)=f_i(x,b)\,, \quad v_t(x,0;b)= 0, \quad i=1,\ldots,l \,.
\end{equation}
The kernel $\widetilde{E} (r,t; 0,b) $ is a diagonal matrix with the
elements $ E_i (r,t; 0,b)$, $i=1,\ldots,l$, which are defined either
by (\ref{E}) with corresponding mass terms \, $m_i$, $i=1,\ldots,l$,
or by (\ref{EM}), in accordance with the value of mass $m_i^2 \geq
n^2/4 $ or $m_i^2 < n^2/4 $, respectively.
Then, the solution  $\Psi $ of the Cauchy problem for the equation
\begin{eqnarray*}
&  &
\Psi_{tt} +   n H \Psi _t - e^{-2Ht} \bigtriangleup \Psi  +  \widetilde{\bf M}\Psi=  0\,
\end{eqnarray*}
with the initial conditions
\[
\Psi (x,0)= \widetilde{\psi} _0 (x),\quad \Psi _t(x,0)=\widetilde{\psi} _1 (x) ,
\]
with the vector-valued functions \, $ \widetilde{\psi}_0 $,  $ \widetilde{\psi}_1 \in C_0^\infty ({\mathbb R}^n) $, $n\geq 2$, can be represented as follows:
\begin{eqnarray*}
\Psi  (x,t)
& = &
e^{-\frac{n-1}{2}t} \widetilde{v}_{\widetilde{\psi }_0}  (x, \phi (t))\\
&  &
+ \,  e^{-\frac{n}{2}t}\int_{ 0}^{1}\big(2  \widetilde{K}_0(\phi (t)s,t)
+ n\widetilde{K}_1(\phi (t)s,t)\big) \widetilde{v}_{\widetilde{\psi }_0}(x, \phi (t)s)\phi (t)\,  ds  \nonumber \\
& &
+\, 2e^{-\frac{n}{2}t}\int_{0}^1
  \widetilde{K}_1(\phi (t)s,t) \widetilde{v}_{\widetilde{\psi } _1 } (x, \phi (t) s)\phi (t)\, ds
, \quad x \in {\mathbb R}^n, \,\, t>0\,,
\end{eqnarray*}
where  $\phi (t):= 1-e^{-t} $ and the kernels $ \widetilde{K}_0 $,
$\widetilde{K}_1 $, are the diagonal matrices with the elements $
\widetilde{K}_{0i}  (z,t)$, $i=1,\ldots,l$, and $ \widetilde{K}_{1i}
(z,t)$, which are defined either by (\ref{K0}) and (\ref{K1})  with
the corresponding mass terms \, $m_i$, $i=1,\ldots,l$, or by the
diagonal matrices with the elements $ \widetilde{K}_{0i}  (z,t;M)$,
$i=1,\ldots,l$, and $ \widetilde{K}_{1i}  (z,t;M)$, which are
defined by (\ref{K0M}) and (\ref{K1M}), in accordance with the value
of mass $m_i^2 \geq n^2/4 $ or $m_i^2 < n^2/4 $, respectively.
 Here, for  the vector-valued function $\widetilde{\psi}  \in C_0^\infty ({\mathbb R}^n)$ and for $x \in {\mathbb R}^n$,
the vector-valued function $\widetilde{v}_{\widetilde{\psi }  }  (x,t)$  is a solution   of the Cauchy problem
$
\widetilde{v}_{tt}-  \bigtriangleup \widetilde{v} =0$, $ \widetilde{v}(x,0)= \widetilde{\psi}  (x)$, $ \widetilde{v}_t(x,0)=0\,.
$
We study the Cauchy problem
    through the integral equation.
To determine that integral equation we  appeal to the operator
\[
\widetilde{G}:= \widetilde{\mathcal K}\circ \widetilde{\mathcal WE} \,,
\]
where the operator $\widetilde{\mathcal WE} $ is defined by (\ref{WEsystem}), that is,
\[
\widetilde{\mathcal WE} [f](x,t;b)= \widetilde v(x,t;b),
\]
and the vector-valued function
$\widetilde v(x,t;b)$
is a solution to the Cauchy problem for the  wave equation,
while $\widetilde {\mathcal K}$ is introduced either by (\ref{1.29l}),
for the large mass matrix, or by (\ref{Pfim}),
for the small mass matrix. Hence,
\[
\widetilde{G}[f]  (x,t)
  =
2   e^{-\frac{n}{2}t}\int_{ 0}^{t} db
  \int_{ 0}^{ e^{-b}- e^{-t}} dr  \,  e^{\frac{n}{2}b}\,\widetilde{E}(r,t; 0,b;M)\widetilde{\mathcal WE} [f](x,r ;b).
\]
Thus, the Cauchy problem (\ref{NWE}), (\ref{ICPHI}) leads to the following integral equation
\begin{eqnarray}
\label{5.1}
\label{4.5}
\Psi (x,t)
 =
\Psi _0(x,t) +
\widetilde{G}[ \widetilde{F}(\Psi ) ] (x,t)    \,.
\end{eqnarray}
Every solution $\Phi =\Phi (x,t) $ to  the  equation (\ref{NWE}) generates the function $\Psi =\Psi (x,t) $, which solves   the last integral equation with
some function $\Psi _0 (x,t)$, that is generated by the solution of the Cauchy problem   (\ref{1.10}).

\paragraph{Solvability of the integral equation  associated with Klein-Gordon equation.}
We are going to apply Banach's fixed-point theorem. In
order to estimate the nonlinear term we use the   Lipschitz
condition (${\mathcal L}$), which
imposes some restrictions on $n$, $\alpha $, $s$. Then we consider
the  equation  (\ref{5.1}), where the vector-valued function
$\Psi _0\in C([0,\infty);L^q ({\mathbb R}^n))$ is given.
The solvability of the integral equation (\ref{5.1}) depends on the
operator $\widetilde G$.
We start with the case of Sobolev space $ H_{(s)} ({\mathbb R}^n)$ with $ s
> n/2$, which is an algebra. In the next theorem the operator
$\widetilde{\mathcal K} $ is generated by the linear part of the
equation (\ref{NWE}).
\begin{thm} \mbox{\rm \cite{birk}}
\label{TIE}
Assume that  $ F(\Psi  )$  is Lipschitz continuous in the  space $H_{(s)} ({\mathbb R}^n)$, $ s > n/2$, and also that $\alpha >0 $.\\
(i) Let the spectrum of the mass matrix $\bf M$ be $\{m_1^2,\ldots
m_l^2\}  \subset (0, (n^2-1)/4] $, and $m=\min\{ m_1, m_2,\ldots ,
m_l \} $. Then  for every given function $ \Psi _0(x ,t) \in
X({R,s,\gamma_0}) $ such that
\begin{eqnarray*}
&  &
\sup_{t \in [0,\infty)}  e^{\gamma_0 t}  \|\Psi _0(x ,t) \|_{H_{(s)} ({\mathbb R}^n)}  < \varepsilon\,, \qquad where \quad
 \gamma_0\leq \frac{n}{2 } -\sqrt{\frac{n^2}{4 }-m^2} \,,
\end{eqnarray*}
and for sufficiently small \,  $\varepsilon $  \, the integral equation (\ref{4.5}) has a unique solution \, $ \Psi  (x ,t) \in X({R,s,\gamma})  $ \,
with $0 < \gamma <    {\gamma_0}/{(\alpha +1)} $. For the solution one has
\[
\sup_{t \in [0,\infty)}  e^{\gamma t}  \|\Psi  (x ,t) \|_{H_{(s)} ({\mathbb R}^n)}  < 2\varepsilon \,.
\]
(ii) If the eigenvalues of the mass matrix are large, $\frac{n}{2 }
\leq  m_i$, $i=1,\ldots,l$, then  for every given function $ \Psi
_0(x ,t) \in X({R,s,0}) $ such that
\[
\sup_{t \in [0,\infty)}    \|\Psi _0(x ,t) \|_{H_{(s)} ({\mathbb R}^n)}  < \varepsilon,
\]
and for sufficiently small   $\varepsilon $   the integral equation (\ref{4.5}) has a unique solution
$ \Psi  (x ,t) $ $\in X({R,s,0})  $, and
\[
\sup_{t \in [0,\infty)}    \|\Psi  (x ,t) \|_{H_{(s)} ({\mathbb R}^n)}  < 2\varepsilon .
\]
\end{thm}

\section{Asymptotic at infinity}
\setcounter{equation}{0}
\renewcommand{\theequation}{\thesection.\arabic{equation}}

For $\varphi  \in C_0^\infty ({\mathbb R}^n)$
let $V_{\varphi }(x, t) $ be a solution     of the Cauchy problem
\[
V_{tt}- \Delta  V =0 , \qquad  V(x,0)= 0 , \quad V_t(x,0)=\varphi  (x).
\]
Denote,
\begin{eqnarray*}
V_{\varphi }^{(k)}(x)
& = &
\frac{(-1)^k}{k!}\left[\left( \frac{\partial }{\partial t} \right)^k V_{\varphi }(x, t) \right]_{t=1}
\in C_0^\infty ({\mathbb R}^n)\,,\quad k=1,2,\ldots\,.
\end{eqnarray*}
Then, for every integer  $N \geq 1$ we have
\begin{eqnarray*}
V_{\varphi }(x,1-e^{-t}) = \sum_{k=0}^{N-1} V_{\varphi }^{(k)}(x)e^{-kt} + R_{V_\varphi ,N}(x,t),  \quad
R_{V_\varphi ,N} \in   C ^\infty \,,
\end{eqnarray*}
where with the constant $ C(\varphi ) $ the remainder  $ R_{V_\varphi ,N} $ satisfies the inequality
\begin{eqnarray*}
| R_{V_\varphi ,N}(x,t) | \leq C(\varphi ) e^{-Nt}\quad  \,\, \mbox{\rm for all}\quad  x \in {\mathbb R}^n \quad  \mbox{\rm and all}\quad  t \in [0,\infty)  \,.
\end{eqnarray*}
Further,  we introduce   the polynomial in $z $ with the smooth in $x \in {\mathbb R}^n$ coefficients  as follows:
\begin{eqnarray*}
\Phi_{asypt}^{(N)}  (x,z)
& = &
 z^{\frac{n-1}{2}}    \sum_{k=0}^{N-1}\left( \frac{n-1}{2}      V_{\varphi_0 }^{(k)}(x)  -(k+1) V_{\varphi_0 }^{(k+1)}(x)
+   V_{\varphi_1 }^{(k)}(x)\right) z^k .
\end{eqnarray*}
where $ x \in {\mathbb R}^n$, $z \in {\mathbb C} $. Thus, we arrive at the next theorem.

\begin{thm} \mbox{\rm \cite{ArXiV-2012}}
\label{Tasymp}
Suppose that  $m= \sqrt{ n^2-1 }/2 $. Then, for every integer positive $N$   the solution of the equation (\ref{25}) with the initial values
$ \varphi_0, \varphi _1  \in C_0^\infty ({\mathbb R}^n) $   has the following asymptotic expansion at infinity:
\begin{eqnarray*}
 \Phi  (x,t) = \Phi_{asypt}^{(N)}  (x,e^{-t})  + O(e^{-Nt-\frac{n-1}{2}t})
\end{eqnarray*}
for large $t$ uniformly for  $x\in {\mathbb R}^n$, in the sense that  for every integer positive $N$  the following
estimate is valid:
\begin{eqnarray*}
\| \Phi  (x,t) - \Phi_{asypt}^{(N)}  (x,e^{-t}) \|_{L^\infty  ({\mathbb R}^n)}
& \leq  &
C(\varphi_0, \varphi _1) e^{-Nt-\frac{n-1}{2}t} \quad \,\, \mbox{  for large }\,\, t\,.
\end{eqnarray*}
\end{thm}
 Unlike to the result by Vasy~\cite{Vasy_2010} the last inequality
does not have the logarithmic term.

\end{document}